\documentclass[12pt]{article}
\usepackage{amsmath, amsthm, amsfonts, makeidx}
\usepackage{multind}

\newlength{\stefan}
\DeclareMathSymbol{\subsetneq}{\mathord}{AMSb}{"26}

\newtheorem{lemma}{Lemma}[section]

\newtheorem{theorem}[lemma]{Theorem}

\newtheorem{proposition}[lemma]{Proposition}
\newtheorem{corollary}[lemma]{Corollary}
\theoremstyle{definition}

\newtheorem{definitie}[lemma]{Definition}
\newtheorem{definition}[lemma]{Definition}
\newtheorem{example}[lemma]{Example}

\renewcommand{\ker}{\mathit{ker}}

\renewenvironment{proof}[1][Proof]{\noindent\textbf{#1.} }{\ \rule{0.5em}{0.5em}}
\input{tcilatex}
\begin{document}

\title{Constructing (almost) rigid rings and a UFD having infinitely
generated Derksen and Makar-Limanov invariant}
\author{David Finston and Stefan Maubach\thanks{%
Funded by Veni-grant of council for the physical sciences, Netherlands
Organisation for scientific research (NWO)}}

\maketitle

\begin{abstract}
An example is given of a UFD which has infinitely generated Derksen
invariant. The ring is \textquotedblleft almost rigid\textquotedblright\
meaning that the Derksen invariant is equal to the Makar-Limanov invariant.
Techniques to show that a ring is (almost) rigid are discussed, among which
is a generalization of Mason's abc-theorem.
\end{abstract}


\section{Introduction and tools}

The Derksen invariant and Makar-Limanov invariant are useful tools to
distinguish nonisomorphic algebras. \ They have been applied extensively in
the context of affine algebraic varieties. Both invariants rely on locally
nilpotent derivations: for $R$ a commutative ring and $A$ a commutative $R$%
-algebra, an $R$-linear mapping $D:A\rightarrow A$ is an $R$-derivation if $%
D $ satisfies the Leibniz rule: $D(ab)=aD(b)+bD(a)$. The derivation $D$ is
locally nilpotent if for each $a\in A$ there is some $n\in \mathbb{N}$ such
that $D^{n}(a)=0$. When $k$ is a field of characteristic $0$ a locally
nilpotent $k$-derivation $D$ of the $k$-algebra $A$ gives rise to an
algebraic action of the additive group of $k,$ $G_{a}(k),$ on $A$ via: 
\begin{equation*}
\exp (tD)(a)\equiv \tsum\limits_{i=0}^{\infty }\frac{t^{i}}{i!}d^{i}(a).
\end{equation*}%
for $t\in k,a\in A.$ Conversely, an algebraic action $\sigma $ of $G_{a}(k)$
on $A$ yields a locally nilpotent derivation $\ $via:%
\begin{equation*}
\frac{\sigma (t,a)-a}{t}|_{t=0}.
\end{equation*}%
In this case, the kernel of $D$ denoted by $A^{D}$ coincides with the ring
of $G_{a}(k)$ invariants in $A.$

The Makar-Limanov invariant of the $R$-algebra $A,$ denoted $ML_{R}(A),$ is
defined as the intersection of the kernels of all locally nilpotent $R$%
-derivations of $A$, while the Derksen invariant, $D_{R}(A)$ is defined as
the smallest algebra containing the kernels of all nonzero locally nilpotent 
$R$-derivations of $A$. The subscript $R$ will be suppressed when it is
clear from the context.

In \cite{R++} the question was posed of whether the Derksen invariant of a
finitely generated algebra over a field could be infinitely generated. In 
\cite{MauPre} an example is given of an infinitely generated Derksen
invariant of a finitely generated $\mathbb{C}$-algebra. In fact, this
example is of a form described in this paper as an \textquotedblleft
almost-rigid ring\textquotedblright : a ring for which the Derksen invariant
is equal to the Makar-Limanov invariant. Despite its simplicity and the
simplicity of the argument, this example has a significant drawback in that
it is not a UFD. In this paper we provide a UFD example having infinitely
generated invariants (it is again an almost-rigid ring).

The paper is organized as follows. Section 1 consists of basic notions and
examples associated with rigidity and almost rigidity. In section 2, the
focus is on rigid and almost rigid rings, with techniques to prove rigidity
or almost rigidity. In section 3, certain rings are shown to be UFDs, and
these are used in section 4 to give the UFD examples having infinitely
generated Makar-Limanov and Derksen invariants.

\textbf{Notations:} If $R$ is a ring, then $R^{[n]}$ denotes the polynomial
ring in $n$ variables over $R$ and $R^{\ast }$ denotes the group of units of 
$R$. The $R$ module of $R$-derivations of an $R$-algebra $A$ is denoted by $%
Der_{R}(A)$ and the set of locally nilpotent $R$-derivations by $LND_{R}(A)$
(the $R$ will be suppressed when it is clear from the context). We will use
the letter $k$ for a field of characteristic zero, and $K$ for an algebraic
closure. The symbol $\partial _{X}\ \ $denotes the derivative with respect
to $X$. When the context is clear, $x,y,z,\ldots $ will represent residue
classes of elements $X,Y,Z,\ldots $ modulo an ideal.

Let $A$ be an $R$-algebra which is an integral domain. Well-known facts that
we need are included in the following:

\begin{lemma}
\label{extra} Let $D\in LND_{R}(A)$. \newline
(1) Then $D(A^{\ast })=0$.\newline
(2) If $D(ab)=0$ where $a,b$ are both nonzero, then $D(a)=D(b)=0$.\newline
(3) If $\ \tilde{D}\in Der_{R}(A)$ and $f\in A$ satisfy $f\tilde{D}$ $\in
LND_{R}(A)$, then $\tilde{D}\in LND_{R}(A)$ and $f\in A^{\tilde{D}}$.
\end{lemma}

\section{(Almost) rigid rings}

As defined in \cite{Fre06} page 196, \cite{CML05}, or \cite{Crach}, a rigid
ring is a ring which has no locally nilpotent derivations except the zero
derivation. \ Examples include the rings {$R:=\mathbb{C}[X,Y,Z]/($}$%
X^{a}+Y^{b}+Z^{c})$ with $a,b,c$ $\geq 2$ and pairwise relatively prime \cite%
{FM06}, and coordinate rings of Platonic $\mathbb{C}^{\ast }$ fiber spaces 
\cite{MM}. We define an almost rigid ring here as a ring whose set of
locally nilpotent derivations is, in some sense, one-dimensional.

\begin{definition}
An $R$-algebra $A$ is called almost-rigid if there is a nonzero $D\in 
\QTR{up}{LND}(A)$ such that $\QTR{up}{LND}(A)=A^{D}D$.
\end{definition}

For a field $F$ any derivation $D$ of $F[X]$ has the form $D=f(X)\partial
_{X}$ . \ Thus the simplest almost-rigid algebra is $F[X]$. Other examples
include the algebras 
\begin{equation*}
{\mathbb{C}[X,Y,Z,U,V]/(}X^{a}+Y^{b}+Z^{c},X^{m}V-Y^{n}U-1)
\end{equation*}%
with $a,b,c$ pairwise relatively prime given in \cite{FM06} as
counterexamples to a cancellation problem. Clearly an almost-rigid algebra
has its Derksen invariant equal to its Makar-Limanov invariant. The
following lemma is useful in determining rigidity.

\begin{lemma}
\label{L1} Let $D$ be a nonzero locally nilpotent derivation on a domain $A$
containing $\mathbb{Q}$. Then $A$\ embeds into $K[S]$ where $K$ is some
algebraically closed field of characteristic zero, in such a way that $%
D=\partial _{S}$ on $K[S]$.
\end{lemma}

\begin{proof}
The proof uses some well-known facts about locally nilpotent derivations.
Since $D\not=0$ is locally nilpotent, we can find an element $p$ such that $%
D^{2}(p)=0,D(p)\not=0$. Set $q:=D(p)$ (and thus $q\in A^{D}$) and observe
that $D$ extends uniquely to a locally nilpotent derivation $\tilde{D}$ of $%
\tilde{A}:=A[q^{-1}]$. Since $\tilde{D}$ has the slice $s:=p/q$ (a slice is
an element $s$ such that $\tilde{D}(s)=1$) we have (see prop.1.3.21 in \cite%
{Essenboek}) $\tilde{A}=\tilde{A}^{\tilde{D}}[s]$ and $\tilde{D}=\partial
_{s}$. Denote by $k\ $\ the quotient field of $\tilde{A}^{\frac{\partial }{%
\partial s}}$ ($=$ quotient field of $A^{D}$) noting that $D$ extends
uniquely to $k[s]$. One can embed $k$ into its algebraic closure $K$, and
the derivation $\partial _{s}$ on $K[s]$, restricted to $A\subseteq K[s]$,
equals $D$.
\end{proof}

As an application, we have

\begin{example}
Let $R:=\mathbb{C}[x,y]=\mathbb{C}[X,Y]/(X^{a}+Y^{b}+1)$ where $a,b\geq 2$.
Then $R$ is rigid.
\end{example}

\begin{proof}
Suppose $D\in \QTR{up}{LND}(R)$, $D\not=0$. Using lemma \ref{L1}, we see $D$
as $\partial _{S}$ on $K[S]\supseteq R$. Now the following lemma
(\textquotedblleft mini-Mason's\textquotedblright ) shows that $x,y$ both
must be constant polynomials in $S$. But that means $D(x)=D(y)=0$, so $D$ is
the zero derivation, contradiction. So the only derivation on $R$ is the
zero derivation, i.e. $R$ is rigid.
\end{proof}

Versions of the following lemma can be found as lemma 9.2 in \cite{Fre06},
and lemma 2 in \cite{lenny}. Here we give it the appellation
\textquotedblleft mini-Mason's\textquotedblright\ as it can be seen as a
very special case of Mason's very useful original theorem. (Note that
Mason's theorem is the case $n=3$ of theorem \ref{Catalan}.)

\begin{lemma}
(Mini-Mason)\label{MiniMason} Let $f,g\in K[S]$ where $K$ is algebraically
closed and of characteristic zero. Suppose that $f^a+g^b\in K^*$ where $%
a,b\geq 2$. Then $f,g\in K$.
\end{lemma}

\begin{proof}
Note that $gcd(f,g)=1$. Taking derivative with respect toNo $S$ gives $%
af^{\prime a-1}=-bg^{\prime b-1}$. So $f$ divides $gg^{\prime }$, so $f$
divides $g^{\prime }$. Same reason, $g$ divides $f^{\prime }$. This can only
be if $f^{\prime }=g^{\prime }=0$.
\end{proof}

Mason's theorem provides a very useful technique in constructing rigid rings
(see \cite{FM06} for an example). With appropriate care, a generalization of
Mason's theorem provides more examples. In this paper, we will use \cite[%
Theorem 2.1]{MichielMason}, which is a corollary of a generalization of
Mason's theorem (see \cite[Theorem 1.5]{MichielMason}).

\begin{theorem}
\label{Catalan} Let $f_{1},f_{2},\ldots ,f_{n}\in K[S]$ where $K$ is an
algebraically closed field containing $\mathbb{Q}$. Assume 
\begin{equation*}
f_{1}^{d_{1}}+f_{2}^{d_{2}}+\ldots +f_{n}^{d_{n}}=0.
\end{equation*}%
Additionally, assume that for every $1\leq i_{1}<i_{2}<\ldots <i_{s}\leq n$, 
\begin{equation*}
f_{i_{1}}^{d_{i_{1}}}+f_{i_{2}}^{d_{i_{2}}}+\ldots
+f_{i_{s}}^{d_{i_{s}}}=0\longrightarrow gcd\{f_{i_{1}},f_{i_{2}},\ldots
,f_{i_{s}}\}=1.
\end{equation*}%
Then 
\begin{equation*}
\sum_{i=1}^{n}\frac{1}{d_{i}}\leq \frac{1}{n-2}
\end{equation*}%
implies that all $f_{i}$ are constant.
\end{theorem}

\begin{example}
Let $R:=\mathbb{C}[X_1,X_2,\ldots,X_n]/(X_1^{d_1}+X_2^{d_2}+\ldots +
X_n^{d_n})$ where $d_1^{-1}+d_2^{-1}+\ldots+d_n^{-1}\leq \frac{1}{n-2}$.
Then $R$ is a rigid ring.
\end{example}

The proof will follow from the more general

\begin{lemma}
\label{kleincatalan} Let $A$ be a finitely generated $\mathbb{Q}$ domain$.$
Consider a subset $\mathcal{F=}\{F_{1},F_{2},\ldots ,F_{m}\}$ of $A$ and
postive integers $d_{1},\ldots d_{n}$ satisfying: 1) $%
P:=F_{1}^{d_{1}}+F_{2}^{d_{2}}+\ldots +F_{m}^{d_{m}}$ is a prime element of $%
A$ and 2) No nontrivial subsum of $F_{1}^{d_{1}},F_{2}^{d_{2}},\ldots
,F_{m}^{d_{m}}$ lies in $(P)$ (e.g. the $F_{i}$ are linearly independent)$.$
\ Additionally, assume that 
\begin{equation*}
d_{1}^{-1}+d_{2}^{-1}+\ldots +d_{n}^{-1}\leq \frac{1}{n-2}.
\end{equation*}%
Set $R:=A/(P)$ and let $D\in LND(R).$ With $f_{i}\in R$ equal to the residue
class of $F_{i}$, we have $D(f_{i})=0$ for all $1\leq i\leq n$.
\end{lemma}

\begin{proof}
Suppose $D\in LND(R)$ where $D\not=0$. Using lemma \ref{L1} with $K$ an
algebraic closure of the quotient field of $R^{D}$ , we realize $D$ as $%
\partial _{S}$ on $K[S]\supseteq R$. In particular, $%
f_{1}(S)^{d_{1}}+f_{2}(S)^{d_{2}}+\ldots +f_{m}(S)^{d_{m}}=0$. By hypothesis
there cannot be a subsum $f_{i_{1}}^{d_{i_{1}}}+f_{i_{2}}^{d_{i_{2}}}+\ldots
+f_{i_{s}}^{d_{i_{s}}}=0$. \ Applying the above theorem \ref{Catalan}, we
find that all $f_{i}$ are constant.
\end{proof}

This lemma also helps in constructing almost-rigid rings not of the form $%
R^{[1]}$ with $R$ rigid.

\begin{example}
\cite{MauPre} Define%
\begin{equation*}
R:=\mathbb{C}[a,b]=\mathbb{C}[A,B]/(A^{3}-B^{2})
\end{equation*}
and 
\begin{equation*}
S:=R[X,Y,Z]/(Z^{2}-a^{2}(aX+bY)^{2}-1).
\end{equation*}
Then $LND(S)=S^{D}D$ where $D:=b\partial _{X}-a\partial _{Y}$.
\end{example}

The following is an example of a rigid unique factorization domain. The
proof of UFD\ property is deferred to the next section.

\begin{example}
\label{example1} Let $n\geq 3$, and in $\mathbb{C}[X_{1},X_{2},\ldots
,X_{n},Y_{1},Y_{2},\ldots ,Y_{n}]$ set 
\begin{equation*}
P:=X_{1}^{d_{1}}+X_{2}^{d_{2}}+\ldots
+X_{n}^{d_{n}}+L_{2}^{e_{2}}+L_{3}^{e_{3}}+\ldots +L_{n}^{e_{n}}
\end{equation*}%
where $L_{i}:=X_{i}Y_{1}-X_{1}Y_{i}$ and 
\begin{equation*}
d_{1}^{-1}+d_{2}^{-1}+\ldots +d_{n}^{-1}+e_{2}^{-1}+e_{3}^{-1}+\ldots
+e_{n}^{-1}\leq 1/(2n-1-2).
\end{equation*}%
Let 
\begin{equation*}
R:=\mathbb{C}[X_{1},X_{2},\ldots ,X_{n},Y_{1},Y_{2},\ldots ,Y_{n}]/(P)
\end{equation*}%
and denote by $x_{i},y_{i,}l_{i}$ the images of $X_{i},Y_{i,}L_{i}$ in $R.$
Then $R$ is an almost-rigid UFD, and $LND(R)=R^{D}D$ where $%
D(x_{i})=0,D(y_{i})=x_{i}.$
\end{example}

\begin{proof}
An elementary argument shows that $R$ is a domain: View 
\begin{equation*}
P\in \mathbb{C}[X_{1},X_{2},\ldots ,X_{n},Y_{1},Y_{2},\ldots
,Y_{n-1}][Y_{n}].
\end{equation*}%
The residue of $P$ modulo $(Y_{1},Y_{2},\ldots ,Y_{n-1})$ has the same
degree in $Y_{n}$ as $P$ and is clearly irreducible. \ 

That any $2n-1$ element subset of $\{x_{i}^{d_{i}},l_{i}^{e_{i}}:1\leq i\leq
n\}$ is algebraically independent over $\mathbb{Q}$ modulo $(P)$ is also
elementary: Suppose that $\sum_{i=1}^{n}X_{i}^{d_{i}}+%
\sum_{i=1}^{n-1}L_{i}^{e_{i}}$ is divisible by $P.$ Lemma \ref{kleincatalan}
yields that for any $E\in LND(R)$ we have $E(x_{i})=0$, and $E(l_{i})=0$. So 
$x_{1}E(y_{i})=x_{i}E(y_{1})$. Since $R$ is a UFD, we can write $%
E(y_{i})=\alpha x_{i}$ for some $\alpha \in R$. So $E=\alpha D$ where $D$ is
as in the statement.
\end{proof}

\section{Factoriality of Brieskorn-Catalan-Fermat rings for $n\geq 5$}

Because of their resemblance to rings arising in Fermat's last theorem, the
Catalan conjecture, and to the coordinate rings of Brieskorn hypersurfaces,
we will call the rings $\mathbb{C}[X_{1},X_{2},\ldots
]/(X_{1}^{d_{1}}+X_{2}^{d_{2}}+\ldots +X_{n}^{d_{n}})$
Brieskorn-Catalan-Fermat (BCF) rings. Our examples depend on the
factoriality of certain BCF rings. \ While the next observation is
undoubtedly well known, a proof is included since we could not find an
explicit one in the literature.

\begin{theorem}
\label{catalanUFD} If $n\geq 5$ and $d_{i}\geq 2$ for all $1\leq i\leq n$,
then $\mathbb{C}[X_{1},X_{2},\ldots ]/(X_{1}^{d_{1}}+X_{2}^{d_{2}}+\ldots
+X_{n}^{d_{n}})$ is a UFD.
\end{theorem}

The result follows from the next two theorems:

\begin{theorem}
\label{Fossum} (Corollary 10.3 of \cite{Fossum}) Let $A=A_0+A_1+\ldots$ be a
graded noetherian Krull domain such that $A_0$ is a field. Let $\mathfrak{m}%
=A_1+A_2+\ldots$. Then $Cl(A)\cong Cl(A_{\mathfrak{m}})$, where $Cl$ is the
class group.
\end{theorem}

\begin{theorem}
\label{local}(\cite{HO74}) A local noetherian ring $(A,\mathfrak{m})$ with
characteristic $A/\mathfrak{m}=0$ and an isolated singularity is a UFD if
its depth is $\geq 3$ and the embedding codimension is $\leq dim(A)-3$.
\end{theorem}

\begin{proof}
(of theorem \ref{catalanUFD}) Write%
\begin{eqnarray*}
A &:&=\mathbb{C}[x_{1},x_{2},\ldots ,x_{n}] \\
&=&\mathbb{C}[X_{1},X_{2},\ldots ]/(X_{1}^{d_{1}}+X_{2}^{d_{2}}+\ldots
+X_{n}^{d_{n}}).
\end{eqnarray*}
Note that by giving appropriate positive weights to the $X_{i}$, the ring $A$
is graded, and $\mathfrak{m}:=A_{1}+A_{2}+\ldots =(x_{1},x_{2},\ldots
,x_{n}),A_{0}=\mathbb{C}$. $A$ now satisfies the requirements of \ref{Fossum}%
, so it is equivalent to show that $A_{\mathfrak{m}}$ is a UFD (note that
\textquotedblleft $A$ ia a UFD\textquotedblright\ is equivalent to
\textquotedblleft $Cl(A)=\{0\}$\textquotedblright ). Now $A_{\mathfrak{m}}$
has only one singularity, namely at the point $\mathfrak{m}$. The ring $A$
is defined by one homogeneous equiation, and therefore, by definition, a
complete intersection. Being a complete intersection implies that the ring $%
A $ is Cohen-Macauley and that its depth is the same as its Krull dimension.
So, the depth of $A$ is $n-1$ which is $\geq 3$ since $n\geq 5$. Now, one
can see $A$ as a subring of the polynomial ring localized at the maximal
ideal $(X_{1},X_{2},\ldots ,X_{n})$. $A$ has codimension 1 in this ring, so
its embedding codimension is 1. $dim(A)-3=n-4$, so, if $n\geq 5$, we have
that the embedding codimension of $A$ equals $1\leq dim(A)-3$. So, if $n\geq
5$, the criteria of \ref{local} are met, and $A_{\mathfrak{m}}$ is a UFD.
\end{proof}

The following lemma of Nagata is a very useful tool in proving factoriality.

\begin{lemma}
\label{nagata} (Nagata) Let $A$ be a domain, and $x\in A$ is prime. If $%
A[x^{-1}]$ is a UFD, then $A$ is a UFD.
\end{lemma}

\begin{lemma}
\label{Rufd} $R$ as in example \ref{example1} is a UFD.
\end{lemma}

\begin{proof}
Note that $X_{2}^{d_{2}}+X_{3}^{d_{3}}+\ldots
+X_{n}^{d_{n}}+(X_{2}Y_{1})^{e_{2}}+(X_{3}Y_{1})^{e_{3}}+\ldots
+(X_{n}Y_{1})^{e_{n}}$ is irreducible for any $d_{i}\geq 1,e_{i}\geq 1$, so $%
R/(x_{1})$ is a domain. Using \ref{nagata} it is enough to show that $%
R[x_{1}^{-1}]$ is a UFD. Define $m_{i}:=y_{i}-\frac{x_{i}}{x_{1}}y_{1}$ for $%
2\leq i\leq n$, and%
\begin{equation*}
S:=\mathbb{C}[x_{1},x_{2},\ldots ,x_{n},m_{2},m_{3},\ldots ,m_{n}].
\end{equation*}
Then $R[x^{-1}]=S[x_{1}^{-1}][Y_{1}]$ where $Y_{1}$ is algebraically
independent over $S[x_{1}^{-1}]$. It is now enough to prove that $S$ is a
UFD. But this follows from theorem \ref{catalanUFD} since $n\geq 3$.
\end{proof}

\section{A UFD having infinitely generated invariants}

\subsection{Definitions}

\begin{definitie}
In $\mathbb{C}^{[7]}=\mathbb{C[}X,Y,Z,S,T,U,V],$ let $%
L_{1}:=Y^{3}S-X^{3}T,L_{2}:=Z^{3}S-X^{3}U,L_{3}:=Y^{2}Z^{2}S-XV$. \ Define $%
P:=X^{d_{1}}+Y^{d_{2}}+Z^{d_{3}}+L_{1}^{d_{4}}+L_{2}^{d_{5}}+L_{3}^{d_{6}}$
where the $d_{i}\geq 2$ are integers . Set 
\begin{equation*}
A:=\mathbb{C}[x,y,z,s,t,u,v]=\mathbb{C}[X,Y,Z,S,T,U,V]/(P),
\end{equation*}%
and let $R$ be the subring $\mathbb{C}[x,y,z].$
\end{definitie}

The elements $s,t,u,v$ in $A$ form a regular sequence; in particular they
are algebraically independent$.$

\begin{definitie}
\begin{equation*}
E:=X^{3}\partial _{S}+Y^{3}\partial _{T}+Z^{3}\partial
_{U}+X^{2}Y^{2}Z^{2}\partial _{V}.
\end{equation*}%
\ Note that $E$ is locally nilpotent and $\ P\in \ker (E).$ Thus $E$ induces
a well defined element of $LND(A)$ denoted by $D$. \ 
\end{definitie}

\subsection{The factoriality of A}

For a 5-tuple of positive integers $\mathbf{d=(}d_{1},d_{2},\ldots ,d_{5})$,
define $Q(\mathbf{d)}%
:=Y^{d_{2}}+Z^{d_{3}}+(Y^{3}S)^{d_{4}}+(Z^{3}S)^{d_{5}}+(Y^{2}Z^{2}S)^{d_{6}} 
$

\begin{proposition}
If $Q(\mathbf{d)}$ is irreducible in $\mathbb{C}[Y,Z,S]$ then $A$ is a UFD.
\end{proposition}

\begin{proof}
Assume that $Q(\mathbf{d)}$ is irreducible. $\ $Note that $A/(x)\cong 
\mathbb{C}[Y,Z,S,T,U,V]/(Q(\mathbf{d)})$ so that $x$ is prime. By Nagata's
lemma \ref{nagata}, it is enough to show that $A[x^{-1}]$ is a UFD. Now
define 
\begin{equation*}
M_{1}:=T-\frac{Y^{3}}{X^{3}}S,M_{2}:=U-\frac{Z^{3}}{X^{3}}S,M_{3}:=V-\frac{%
Y^{2}Z^{2}}{X}S,
\end{equation*}%
write $m_{i}$ for the image of $M_{i}$ in $A[x^{-1}]$, and let 
\begin{equation*}
B=\mathbb{C}[x,y,z,m_{1},m_{2},m_{3}][x^{-1}].
\end{equation*}%
Since $D(s)$ $=x^{3},$ $\frac{s}{x^{3}}$ is a slice for the extension of $D$
to $A[x^{-1}]=B[s],$ with $s$ transcendental over $B$. Consider $C:=\mathbb{C%
}%
[X,Y,Z,M_{1},M_{2},M_{3}]/(X^{d_{1}}+Y^{d_{2}}+Z^{d_{3}}+M_{1}^{d_{4}}+M_{2}^{d_{5}}+M_{3}^{d_{6}}) 
$. This ring is a UFD by theorem \ref{catalanUFD}, so $C[x^{-1}]=B$ is also
a UFD, \ from which we deduce that $B[s]=A[x^{-1}]$ is a UFD.
\end{proof}

The polynomial $Q(\mathbf{d})$ is irreducible for infinitely many positive
integer choices of the $d_{i}$; take for example $\gcd (d_{2},d_{3})=1$ and $%
d_{2}\geq \max (3d_{4},2d_{6})$.

\subsection{$A$ is not finitely generated}

In this section, we assume that $d_{1},\ldots ,d_{6}$ are such that $Q(%
\mathbf{d})$ is irreducible (i.e. $A$ is a UFD), and such that $%
d_{1}+d_{2}+\ldots +d_{6}\leq \frac{1}{4}$ (note that by neccessity $%
d_{1,}d_{2},d_{3}\geq 4)$. The following lemma shows that $A$ is an
almost-rigid ring.

\begin{lemma}
\label{L2}Any locally nilpotent derivation on $A$ is a multiple of $D$.
\end{lemma}

\begin{proof}
Let $\triangle $ be a nonzero LND on $A$. By lemma \ref{kleincatalan}, since
we assumed $\sum_{i=1}^{6}d_{i}\leq \frac{1}{4}$, we see that $%
x,y,z,l_{1},l_{2},l_{3}$ must be in $A^{\triangle }$. So $\triangle
(l_{1})=0 $, so $x^{3}\triangle (t)=y^{3}\triangle (s)$, and thus $\triangle
(S)=x^{3}\alpha $ for some $\alpha \in A$ (since $A$ is a UFD). Using $%
\triangle (l_{1})=\triangle (l_{2})=\triangle (l_{3})=0$ this yields $%
\triangle (T)=y^{3}\alpha ,\triangle (U)=z^{3}\alpha ,\triangle
(V)=x^{2}y^{2}z^{2}\alpha $, i.e.$\triangle =\alpha D$.
\end{proof}

\begin{lemma}
\label{L5} $A^D\subseteq (x,y,z)A+R$.
\end{lemma}

\begin{proof}
Let 
\begin{eqnarray*}
\mathcal{J} &:&=(X^{3},Y^{3},Z^{3},X^{2}Y^{2}Z^{2})(X,Y,Z)\mathbb{C}^{[7]},
\\
H &:&=(x,y,z)A\supseteq J:=(x^{3},y^{3},z^{3},x^{2}y^{2}z^{2})H.
\end{eqnarray*}
Both $J$ and $H$ are $D$ stable ideals of $A$. Denote by $\bar{D}$ the
locally nilpotent derivation induced by $D$ on $\overline{A}:=A/J,\overline{H%
}:=H/J,$ and $\overline{R}\ $the image of $R$ in $\overline{A}$. Note that $%
\bar{D}(\overline{H})=0$. We will prove that $\bar{A}^{\bar{D}}\subseteq 
\bar{H}+\bar{R}$, which will imply that $A^{D}+J\subseteq H+J+R$, and the
required result then follows since $J\subseteq H$.

To that end assume there exists $h\in \overline{A}^{\overline{D}}$ with $%
h\not\in \bar{H}+\bar{R}$. Note that since $P\in \mathcal{J}$ we have 
\begin{equation*}
\overline{A}\cong (\mathbb{C}^{[7]}/(P))/(\mathcal{J}/(P))\cong \mathbb{C}%
^{[7]}/\mathcal{J}.
\end{equation*}
With $\overline{x},\overline{y},\overline{z},\overline{s},\overline{t},%
\overline{u},\overline{v}$ denoting as usual the images of $S,T,U,V$ in $%
\overline{A,}$ we have $\overline{A}=\bar{R}[\overline{s},\overline{t},%
\overline{u},\overline{v}]$, a polynomial ring over $\overline{R}.$

Assign degree $0$ to elements of $\overline{R}$, weights $wt(\overline{v}%
)>>wt(\overline{u)}>>wt(\overline{t})>>wt(\overline{s}),$ and well order
monomials $\overline{s}^{a}\overline{t}^{b}\overline{u}^{c}\overline{v}^{d}$
in $\overline{A}$ lexicographically. \ \ By assumption there exists a
monomial $M$ of lowest order appearing in $h$ which is not in $\bar{H}+\bar{R%
}$ . Say $M:=\overline{r}\overline{s}^{a}\overline{t}^{b}\overline{u}^{c}%
\overline{v}^{d}$ where $r\in \bar{R}\backslash \bar{H}$.

First assume $d\not=0$. Since $\bar{D}(h)=0$, the nonzero monomial $d%
\overline{x}^{2}\overline{y}^{2}\overline{z}^{2}r\overline{s}^{a}\overline{t}%
^{b}\overline{u}^{c}\overline{v}^{d-1}$ must appear in the $\bar{D}$%
-derivative of at least one other monomial $N$ occurring in $h$. Notice that
then $N$ must also have $\bar{R}$-coefficient not in $\bar{H}$, as otherwise 
$\bar{D}(N)=0$ (since $\bar{D}(\bar{H})=0$). Since $\bar{D}N$ contains the
monomial $\ d\overline{x}^{2}\overline{y}^{2}\overline{z}^{2}r\overline{s}%
^{a}\overline{t}^{b}\overline{u}^{c}\overline{v}^{d-1}$, $\bar{D}N$ has
degree $a+b+c+d-1$. But the derivation $\bar{D}$ decreases degree by exactly
one, so that $N$ must have degree $a+b+c+d$. Since $M$ was the lowest degree
polynomial with lowest possible lexicographic ordering, $N$ then must have a
higher lexicographic ordering than $\overline{s}^{a}\overline{t}^{b}%
\overline{u}^{c}\overline{v}^{d}$. But then all (four) terms in $\bar{D}(N)$
will have higher lexicographic ordering than $\overline{s}^{a}\overline{t}%
^{b}\overline{u}^{c}\overline{v}^{d-1}\ $. So, such a monomial $N$ will not
exist, which is a contradiction for this case.

The cases where $d=0,c\not=0$, and $d=c=0,b\not=0$, and $d=c=b=0,a\not=0$ go
similarly, leading to a contradiction. ($d=c=b=a=0$ implies $M\in \bar{R}$,
which we excluded). So, the assumption that $h\not\in \bar{H}+\bar{R}$, was
wrong. Thus $h\in \bar{H}+\bar{R}$ as claimed.
\end{proof}

\begin{lemma}
\label{L6} For each $n\in \mathbb{N}$, there exists $F_{n}\in A^{D}$ which
satisfies $F_{n}=xV^{n}+f_{n}$ where $f_{n}\in
\sum_{i=0}^{n-1}R[s,t,u]v^{i}\subset A$.
\end{lemma}

\begin{proof}
It is shown in several places, for example \cite{scriptieIk}, \cite{DF}, or
page 231 of \cite{Essenboek}, that already on $\mathbb{C}^{[7]}$ there exist
such $\tilde{F}_{n}$ which are in the kernel of the derivation $E$ (they are
key to the proof that the kernel of $E$ is not finitely generated as a $%
\mathbb{C}$-algebra, and therefore yields a counterexample to Hilbert's 14th
problem). \ By taking for $F_{n}$ the image of $\tilde{F}_{n}\ $in $A$ we
obtain the desired kernel elements.
\end{proof}

\begin{corollary}
\label{C2} $A^D$ is not finitely generated as a $\mathbb{C}$-algebra.
\end{corollary}

\begin{proof}
Suppose $A^{D}=R[g_{1},\ldots ,g_{s}]$ for some $g_{i}\in A$. Since $%
A^{D}\subseteq R+(x,y,z)$ by lemma \ref{L5}, we can assume that all $%
g_{i}\in (x,y,z)$. Define $\mathcal{F}_{n}(A):=\sum_{i=0}^{n-1}R[S,T,U]V^{i}$
which is a subset of $A$. Choose $n$ such that $g_{i}\in \mathcal{F}_{n}(A)$
for all $1\leq i\leq s$. Now $F_{n}\in \mathcal{F}_{n}(A)\cap A^{D}$. Then $%
F_{n}=P(g_{1},\ldots ,g_{s})$ for some $P\in R^{[s]}$. Compute modulo $%
(x,y,z)^{2}$. Since each $g_{i}\in (x,y,z)$, we have 
\begin{equation*}
P(g_{1},\ldots ,g_{n})\equiv r_{1}g_{1}+\ldots +r_{n}g_{n}\func{mod}%
(x,y,z)^{2}
\end{equation*}
for some $r_{i}\in R$. So $F_{n}\in Rg_{1}+\ldots +Rg_{n}+(x,y,z)^{2}$. In
particular, $F_{n}\in \mathcal{F}_{n}(A)+(x,y,z)^{2}$. Notice that $%
F_{n}-xV^{n}\in \mathcal{F}_{n}(A)\subseteq \mathcal{F}_{n}(A)+(x,y,z)^{2}$,
so that $xV^{n}\in \mathcal{F}_{n}(A)+(x,y,z)^{2}$. But this is obviously
not the case, contradicting the the assumption that \textquotedblleft $%
A^{D}=R[g_{1},\ldots ,g_{s}]$ for some $g_{i}\not\in R$\textquotedblright .
Thus $A^{D}$ is not finitely generated as an $R$-algebra, a fortiori as a $%
\mathbb{C}$-algebra.
\end{proof}

Using lemma \ref{L2} we know that there is only one kernel of a nontrivial
LND on $A$, so the following result is obvious.

\begin{corollary}
\label{C3} $ML(A)=Der(A)=A^D$ is not finitely generated.
\end{corollary}

{\em 
\begin{tabular}{ll}
Department of mathematics & Department of Mathematics\\
New Mexico State University & Radboud University Nijmegen\\
Las Cruces,  USA & Nijmegen, The Netherlands\\
email: dfinston@nmsu.edu & email: s.maubach@math.ru.nl\\
\end{tabular}
}

\end{document}